\documentclass[12pt]{article}
\usepackage{amssymb,latexsym}
\usepackage{amsmath}
\def\RR{\mathbb{R}}

\begin{document}


\title{\bf Univariate spline quasi-interpolants and applications to
numerical analysis}

\author{ Paul Sablonni\`ere\\
INSA and IRMAR, Rennes}
\date{}
\maketitle

\begin{abstract}
We describe some new univariate spline quasi-interpolants on uniform
partitions of bounded
intervals. Then we give some  applications to numerical analysis:
integration, differentiation and
approximation of zeros.
\end{abstract}
\section{Introduction}

Univariate spline {\sl quasi-interpolants} (abbr. QIs)  can be defined as
operators of the form
$$
Qf=\sum_{j\in J} \mu_j(f) B_j
$$
where $\{B_j, j\in J\}$ is the B-spline basis of some space of splines, say
of degree $d$, on
a bounded interval $I=[a,b]$ endowed with some partition $X_n$ of $I$ in $n$
subintervals. We denote by
$\Pi_d$ the space of polynomials of total degree at most $d$. In general we
impose that $Q$ is
exact on the space $\Pi_d$, i.e. $Qp=p$ for all $p\in \Pi_d$. Some authors
impose further that $Q$
is a projector on the space of splines itself (see e.g.\cite{dvl},
\cite{lls},\cite{ls}). As a
consequence of this property, the approximation order is $O(h^{d+1})$ on
smooth functions, $h$
being the maximum steplength of the partition $X_n$.  The coefficients
$\mu_j$ are local linear functionals which are in general of one of the
following types:\\

(i) {\sl differential type} : $\mu_j(f)$ is a linear combination of values
of derivatives of $f$, of
order at most $d$, at some point in $supp(B_j)$ (see e.g. \cite{db1},
\cite{db2}). The associated
quasi-interpolant is called a {\sl differential quasi-interpolant} (abbr.
DQI).\\

(ii) {\sl integral type} : $\mu_j(f)$ is a linear combination of weighted
mean values of $f$,
i.e. of quantities $\int_a^b fw_j$ where $w_j$ can be, for example, a
linear combination of
B-splines (see e.g. \cite{db1}, \cite{db2}, \cite{ssb}). The associated
quasi-interpolant is called an {\sl integral quasi-interpolant} (abbr. iQI).\\

(iii) {\sl discrete type} : $\mu_j(f)$ is a linear combination of discrete
values of $f$ at some
points in the neighbourhood of $supp(B_j)$ (see e.g. \cite{ccl},
\cite{lls}, \cite{ls}, \cite{ps1}).
The associated quasi-interpolant is called a {\sl discrete
quasi-interpolant} (abbr. dQI).\\

The main advantage of QIs is that they have a direct construction without
solving any system of
linear equations. Moreover,they are local, in the sense that the value of
$Qf(x)$ depends only on
values of $f$ in a neighbourhood of $x$. Finally, they have a rather small
infinity norm, so they
are nearly optimal  approximants.\\

In this paper, we only consider dQs, neither DQIs nor iQIs. We also
restrict our study to
splines defined on {\sl uniform partitions} of $I=[a,b]$. Our aim is to
give {\sl explicit formulas}
for dQIs of degrees $2\le d\le 5$ and some applications to three classical
problems in numerical
analysis: approximate integration and derivation and location of zeros of
functions. The paper is
organised as follows. In sections 2, we recall some facts about splines and
quasi-interpolants. In section 3, we describe dQIs of degrees $2\le d\le 5$
and we give their infinity norms and their approximation orders. In
sections 4 and 5, we give the
associated {\sl quadrature formulas} and some numerical examples. In
section 6, we give the
{\sl differentiation matrices} for quadratic and cubic splines with
numerical examples. Finally, in
section 7, we show, on a simple polynomial example, how quadratic dQIs can
be applied to the
{\sl location of zeros} of functions.


\section{ Spline spaces on uniform  partitions and dQIs}

For $I=[a,b]$, we denote by  $S_d(I,X_n)$ the space of splines of degree
$d$ and class $C^{d-1}$ on the uniform partition $X_n=\{x_i=a+ih,\; 0\le
i\le n\}$ with meshlength
$h=\frac{b-a}{n}$. A basis of this space is  $\{B_j, j\in J\}$,
with $J=\{1,2,\ldots,n+d\}$. With these notations,
$supp\,(B_j)=[x_{j-d-1},x_j]$,
and ${\cal{N}}_j=\{x_{j-d},\ldots,x_{j-1}\}$ is the set of the $d$ interior
knots in
the support of $B_j$.
As usual, we add multiple knots at the endpoints: $a=x_0=x_{-1}=\ldots=x_{-d}$
and $b=x_n=x_{n+1}=\ldots=x_{n+d}$.
We recall the represention of monomials in terms of symmetric functions of knots
in ${\cal{N}}_j$ \cite{dvl},\cite{sc}.
$$
e_r(x)=x^r=\sum_{j\in J}\theta_j^{(r)}B_j(x),\quad \theta_j^{(r)}={d\choose
r}^{-1}\, symm_r({\cal{N}}_j),\;\; 0\le r \le d.
$$
In particular, the {\sl Greville points}
$$
\theta_j=\theta_j^{(1)}=\frac{1}{d}\sum_{\ell=1}^{d} x_{j-\ell}
$$
are the coefficients of $e_1=\sum_{j\in J}\theta_j B_j$
and the vertices of the {\sl control polygon} of $S=\sum_{j\in
J}c_j B_j$ are the {\sl control points}
$\mathcal{P}=\{\tilde c_j=(\theta_j, c_j), j\in J\}$.
The  {\sl Schoenberg-Marsden} operator is the simplest  discrete
quasi-interpolant
which is exact on the space $\Pi_1$:
$$
S_1f=\sum_{j\in J}f(\theta_j) \,B_j,\;\;\;\;
S_1p=p\quad \forall p\in \Pi_1.
$$

----------------

A {\sl discrete quasi-interpolant} (abbr. dQI) of degree $d$
is a spline operator of the form:
$$
Q_df=\sum_{j\in J}\mu_j(f)\,B_j
$$
whose coefficients $\mu_j(f)$ are linear combinations
of values of $f$ on either the set \\
 $T_n$ ( for $d$ even) or on the set $X_n$  (for $d$ odd), where
$$
T_n=\{t_j, j\in J\},\;\; t_j=\frac12(x_{j-2}+x_{j-1}),\quad
X_n=\{x_j,\;\; 0\le j\le n\}.
$$
Therefore, for $d$ even, we set $f(T_n)=\{f_j=f(t_j),\;\; j\in J\}$,
and for  $d$ odd, we set
$f(X_n)=\{f_j=f(x_j),\;\; 0\le j\le n \}$.
Moreover, $Q_d$ {\sl is exact on} $\Pi_d$:
$$
Q_d\,p=p \qquad\forall p\in \Pi_d.
$$
For the construction of dQIs, i.e. for the determination of functionals
$\{\mu_j(f), j\in J\}$, the exactness of $Q_d$ on $\Pi_d$ amounts to solve
a system of linear
equations for interior B-splines and a finite number of specific linear
systems for boundary
B-splines. The determinants of these systems being Vandermonde
determinants, there is existence and
unicity of dQIs with the above assumptions (see also \cite{ccl} for more
general cases). For the
sake of completeness, we give below complete formulas for degrees $2\le
d\le 5$. Moreover, we give
exact values or upper bounds of $\Vert Q_d \Vert_{\infty}$ and their
approximation order on smooth
functions. Actually, it is well known (see e.g. \cite{dvl}, chapter 5) that
for any subinterval
$I_k=[x_{k-1},x_k],\, 1\le k\le n$, and for any function $f$
$$
\Vert f-Q_df \Vert_{\infty,I_k}\le (1+\Vert Q_d \Vert_{\infty})
d_{\infty,I_k}(f,\Pi_d)
$$
where the distance of $f$ to polynomials is defined by
$$
d_{\infty,I_k}(f,\Pi_d)=\inf\{\Vert f-p \Vert_{\infty,I_k}, p\in\Pi_d\}
$$
Here, as usual, $\Vert f-p \Vert_{\infty,I_k}=max_{x\in I_k} \vert
f(x)-p(x) \vert$.
Therefore, for $f$ smooth enough, e.g. $f\in C^{d+1}(I)$, this implies that
$\Vert f-Q_df \Vert_{\infty}=O(h^{d+1})$.

\section{ Discrete Quasi-Interpolants of degrees $2\le d\le 5$}
\subsection{$C^1$ Quadratic dQI}

For the $C^1$-quadratic dQI $Q_2f=\sum_{j=1}^{n+2}\mu_j(f) \,B_j$,
the coefficient functionals are easy to compute (details are given in
\cite{ps2},\cite{ps3}):\\

$\mu_1(f)=f_1$,
$\mu_2(f)=\frac16(-2f_1+9f_2-f_3)$,
$\mu_{n+1}(f)=\frac16(-f_n+9f_{n+1}-2f_{n+2})$,\\
$\mu_{n+2}(f)=f_{n+1}$, and for $3\le j\le n$
$$
\mu_j(f)=\frac18(-f_{j-1}+10 f_j-f_{j+1})
$$
The exact value $\Vert Q_2 \Vert_{\infty}=1.4734$ has been computed in
\cite{ps3}.
Therefore, for $f\in C^3(I)$ for example, we have the following error estimates
$$
\Vert f-Q_2f \Vert_{\infty,I_k}\le \frac52 d_{\infty,I_k} (f,\Pi_2)\; \mbox{for}\;
 {1\le k\le n}\Longrightarrow \Vert f-Q_2f \Vert_{\infty}=O(h^3).
$$


\subsection{$C^2$ Cubic dQI }

For the $C^2$ cubic dQI $Q_3f=\sum_{j=1}^{n+3}\mu_j(f) \,B_j$,
the coefficient functionals are respectively:

$\mu_2(f)=\frac{1}{18}(7f_0+18f_1-9f_2+2f_3)$,
$\mu_{n+2}(f)=\frac{1}{18}(2f_{n-3}-9f_{n-2}+18f_{n-1}+7f_n)$,\\
$\mu_1(f)=f_0,\quad \mu_{n+3}(f)=f_{n}$, and for $3\le j\le n+1$
$$
\mu_j(f)=\frac16(-f_{j-3}+8 f_{j-2}-f_{j-1}).
$$
As $\vert\mu_2\vert_{\infty}=\vert\mu_{n+2}\vert_{\infty}=2$ and
 $\vert\mu_j\vert_{\infty}= \frac53$ for $3\le j\le n+1$,
we obtain the upper bound $\Vert Q_3 \Vert_{\infty}\le 2$.
It is possible to improve that result by writing the operator in the
"quasi-Lagrange" form:
$$
Q_3f=\sum_{j=1}^{n+3}f_j L_j
$$
where the {\sl fundamental functions} are linear combinations of B-splines,
e.g. for $4\le j\le n$, $L_j=\frac16(-B_{j+3}+8B_{j+2}-B_{j+1})$.
It is well known that $\Vert Q_3\Vert_{\infty}$ is equal
to the Chebyshev norm of the Lebesgue function:
$$
\Lambda_3=\sum_{j=1}^{n+3} \vert L_j \vert.
$$
In each interval of the uniform partition, $\Lambda_3$ is bounded above by
the cubic polynomial
whose Bernstein-B\'ezier (abbr. BB-) coefficients are sums of absolute
values of BB-coefficients of
fundamental functions. This allows to see that the maximum of
$\Lambda_3$ is attained  in the interval $[x_1,x_2]$ and we obtain:
$$
\Vert Q_3 \Vert_{\infty}=\Vert \Lambda_3 \Vert_{\infty}\approx 1.631.
$$
From that we deduce for $f\in C^4(I)$ for example, we have the following
error estimates
$$
\Vert f-Q_3f \Vert_{\infty,I_k}\le \frac83 d_{\infty,I_k}(f,\Pi_3)\;
\mbox{for}\; 
{1\le k\le n}\Longrightarrow \Vert f-Q_3f \Vert_{\infty}=O(h^4).
$$

\subsection{ Quartic dQI}

For the $C^3$ quartic dQI
$Q_4f=\sum_{j=1}^{n+4}\mu_j(f) \,B_j$,
the coefficient functionals are respectively:\\

 $\mu_1(f)=f_1,\quad\mu_{n+4}(f)=f_{n+2}$,\\
$\mu_2(f)=\frac{17}{105}f_1+\frac{35}{32}f_2
-\frac{35}{96}f_3
+\frac{21}{160}f_4-\frac{5}{224}f_5$,\\
$\mu_3(f)=-\frac{19}{45}f_1+\frac{377}{288}f_2
+\frac{61}{288}f_3-\frac{59}{480}f_4+\frac{7}{288}f_5$,\\
$\mu_4(f)=\frac{47}{315}f_1-\frac{77}{144}f_2+\frac{251}{144}f_3
-\frac{97}{240}f_4+\frac{47}{1008}f_5$,\\
$\mu_{n+1}(f)=\frac{47}{315}f_{n+2}-\frac{77}{144}f_{n+1}+\frac{251}{144}f_n
-\frac{97}{240}f_{n-1}+\frac{47}{1008}f_{n-2}$,\\
$\mu_{n+2}(f)=-\frac{19}{45}f_{n+2}+\frac{377}{288}f_{n+1}
+\frac{61}{288}f_n-\frac{59}{480}f_{n-1}+\frac{7}{288}f_{n-2}$,\\
$\mu_{n+3}(f)=\frac{17}{105}f_{n+2}+\frac{35}{32}f_{n+1}-\frac{35}{96}f_n
+\frac{21}{160}f_{n-1}-\frac{5}{224}f_{n-2}$,\\

and for $5\le j\le n$
$$
\mu_j(f)=\frac{47}{1152}(f_{j-4}+f_j)-
\frac{107}{288}(f_{j-3}+f_{j-1})+\frac{319}{192}f_{j-2}
$$
Let us give some details on the computation of functionals
$\mu_k, k=2,3,4$ . As $\mu_k(e_r)=\theta_k^{(r)}$ for $0\le r\le 4$, we
determine the five
coefficients of the discrete functional
$$
\mu_k(f)=\alpha_kf_1+\beta_kf_2+\gamma_kf_3+\delta_kf_4+\zeta_kf_5
$$
as solutions of the three corresponding linear systems ($2\le k\le 4$) of
$5\times 5$ linear
equations
$$
t_1^r\alpha_k+t_2^r\beta_k+t_3^r\gamma_k+t_4^r\delta_k+t_5^r\zeta_k=\theta_k
^{(r)},\;\; 0\le r\le 4
$$
They have the same Vandermonde determinant $V_5(t_1,t_2,t_3,t_4,t_5)\neq 0$
since the $t_i's$
are distinct. Therefore they have unique solutions. The same technique is
applied to the
computation of other coefficient functionals.\\

As $\vert \mu_2 \vert_{\infty}=\vert \mu_{n+3}\vert_{\infty}\approx 1.77$,
$\vert \mu_3 \vert_{\infty}=\vert \mu_{n+2} \vert_{\infty}\approx 2.09$,
$\vert \mu_4 \vert_{\infty}=\vert \mu_{n+1} \vert_{\infty}\approx 2.88$,
$\vert \mu_j \vert_{\infty}\approx 2.49$ for $1\le j\le 5$,  we can
conclude that
\medskip
$\Vert Q_4 \Vert_{\infty}\le 2.88$
and that for $f\in C^5(I)$ for example, we have the following error estimates
$$
\Vert f-Q_4f \Vert_{\infty,I_k}\le 4\, d_{\infty,I_k}(f,\Pi_4)\; \mbox{for}\; 
{1\le k\le n}\Longrightarrow \Vert f-Q_4f \Vert_{\infty}=O(h^5).
$$

\subsection{Quintic dQI}

For the $C^4$ quintic dQI $Q_5f=\sum_{j=1}^{n+5}\mu_j(f) \,B_j$,
the coefficient functionals are respectively:
\medskip

$\mu_1=f_0,\;\mu_{n+5}=f_n$,\\
$\mu_2=\frac{163}{300}f_0+f_1-f_2+\frac23
f_3-\frac14 f_4+\frac{1}{25} f_5$\\
$\mu_3=\frac{1}{200}f_0+\frac{103}{60}f_1
-\frac{73}{60}f_2+\frac{7}{10}f_3-\frac{29}{120}
f_4+\frac{11}{300} f_5$\\
$\mu_4=-\frac{41}{400}f_0+\frac{43}{60}f_1
+\frac{103}{120}f_2-\frac{7}{10}f_3+\frac{13}{48}
f_4-\frac{13}{300} f_5$\\
( symmetric formulas for $n+2\le j\le n+4$),
and for $5\le j\le n$:
$$\mu_j=
\frac{13}{240}(f_{j-5}+f_{j-1})-\frac{7}{15}(f_{j-4}+f_{j-2})+\frac{73}{40}f
_{j-3}
$$
As $\vert \mu_2 \vert_{\infty}=\vert \mu_{n+4} \vert_{\infty}=3.5$,
 $\vert \mu_3 \vert_{\infty}=\vert \mu_{n+3} \vert_{\infty}\approx 3.92$,
 $\vert \mu_4 \vert_{\infty}=\vert \mu_{n+2} \vert_{\infty}\approx 2.69$, and
 $\vert \mu_j \vert_{\infty}\approx 2.87$ for  $5\le j\le n$, we deduce that
$\Vert Q_5\Vert_{\infty}\le 3.92.$ Using a similar technique as for cubics,
we find that
$\Vert Q_5\Vert_{\infty}\approx 3.106.$
Therefore, for $f\in C^6(I)$ for example, we have the following error estimates
$$
\Vert f-Q_5f \Vert_{\infty,I_k}\le 4.5\, d_{\infty,I_k}(f,\Pi_5)\;
\mbox{for}\; 
{1\le k\le n}\Longrightarrow \Vert f-Q_5f \Vert_{\infty}=O(h^6).
$$

\section{Application to numerical integration}

Newton-Cotes formulas are obtained by integrating interpolation polynomials
(see e.g.
\cite{dr},\cite{e},\cite{ku}). In the same way, integrating spline
quasi-interpolants give
interesting quadrature formulas (abbr. QF)  which are easily deduced from
the above computations.
We use the notations
$$
\mathcal{I}(f)=\int_a^b f,\;\; \mathcal{I}_d(f)=\int_a^b Q_df=\sum_{j\in J}
\mu_j(f)\int_a^b B_j,
\;\; E_d(f)=\mathcal{I}(f)-\mathcal{I}_d(f).
$$
As $\int_a^b B_j=\frac{1}{d}(x_{j}-x_{j-d-1})$ and $\mu_j(f)$ are known
explicitly, we can compute
the following quadrature formulas. Moreover, as QIs give the best
approximation order, we can
conclude that $ E_d(f)=O(h^{d+1})$ for $f\in C^{d+1}(I)$, where $h$ is the
meshlength.
Moreover, as for Newton-Cotes formulas, we get a higher approximation order
for even degrees.\\

{\bf (i) QF for quadratics}\\

$\mathcal{I}_2(f)=\int_a^b Q_2f=h\sum_{j=4}^{n-1}
f_j+h\left[\frac19(f_1+f_{n+2})+
\frac78(f_2+f_{n+1})+\frac{73}{72}(f_3+f_{n})\right]$\\

Error : for $f\in C^4(I)$,
$E_2(f)=\mathcal{I}(f)-\mathcal{I}_2(f)$
$=\frac{23}{5760}\,h^4 D^4f(c)-\frac{1}{192}\,h^5 D^4f(\bar c)$\\
This result is proved in \cite{ps3}.\\

Error for Simpson:
$E_2^*(f)={\mathcal{I}}(f)-{\mathcal{I}}_2^*(f)=-\frac{1}{180}\,h^4 D^4f(c)$\\

By comparing the two above errors, we see that the linear combination
(extrapolation):
$$\tilde{\mathcal{I}}_2(f)=
\frac{1}{55}(32{\mathcal{I}}_2(f)+23{\mathcal{I}}_2^*(f))
$$
is such that ${\mathcal{I}}(f)-\tilde{\mathcal{I}}_2(f)=O(h^5)$.\\

{\bf  (ii) QF for cubics}\\

${\mathcal{I}}_3(f)=\int_a^b Q_3 f=h\sum_{j=4}^{n-4}
f_j+h\left[\frac{23}{72}(f_0+f_{n})+\frac{4}{3}(f_1+f_{n-1})+
\frac{19}{24}(f_2+f_{n-2})+\frac{19}{18}(f_3+f_{n-3})\right]$\\

 Error:
$E_3(f)={\mathcal{I}}(f)-{\mathcal{I}}_3(f)=O(h^4)$
for $f\in C^4(I)$. Numerical experiments show that this formula is not as
good as the preceding
one.\\

{\bf (iii) QF for quartics}\\

 $\int_a^b Q_4 f={\mathcal{I}}_4(f)=h\sum_{j=6}^{n-3}
f_j+h\left[\frac{206}{1575}(f_1+f_{n+2})
+\frac{107}{128}(f_2+f_{n+1})+\frac{6019}{5760}(f_3+f_{n})\right]$\\
$+h\left[\frac{9467}{9600}(f_4+f_{n-1})+
\frac{13469}{13440}(f_5+f_{n-2})\right]$\\

 Error:
${\mathcal{I}}(f)-{\mathcal{I}}_4(f)=O(h^6)$ for  $f\in C^6(I)$. This is a
remarkable formula, which
can be compared to the Newton-Cotes formula of the same order. Numerical
experiments show that the
error for the former QF has also the opposite sign of the error for the
latter, as in the quadratic
case. The proof will be given elsewhere.\\

{\bf (iv) QF  for quintics}\\

$\mathcal{I}_5(f)=\int_a^b Q_5f=h\sum_{j=6}^{n-6}
f_j+h\left[\frac{157}{480}(f_0+f_{n})+
\frac{961}{720}(f_1+f_{n-1})+
\frac{133}{180}(f_2+f_{n-2})\right]+$\\
$h\left[\frac{271}{240}(f_3+f_{n-3})
+\frac{1393}{1440}(f_4+f_{n-4})+
\frac{361}{360}(f_5+f_{n-5})\right]$\\

 Error:
$E_5(f)=\mathcal{I}_(f)-\mathcal{I}_5(f)=O(h^6)$
for $f\in C^6(I)$. Numerical experiments show that this formula is not as
good as the preceding
one.\\

\section{Numerical examples}

We compare numerical results on QF applied to the computation of
$$
\mathcal{I}(f_1)=\int_{-1}^1 \frac{1}{1+16x^2}\,dx\;\; \mbox{and} \;\;
\mathcal{I}(f_2)=\int_{-1}^1 e^{-x}sin(5\pi x)\,dx.
$$

{\bf (i) QF/dQI degrees 2 and 3}\\

$E_2(f)=\mathcal{I}(f)-\mathcal{I}_2(f)=O(h^4)$,
$E_3(f)=\mathcal{I}(f)-\mathcal{I}_3(f)=O(h^4)$
for $f\in C^4(I)$\\

{\bf Simpson's QF}
$E_2^*(f)=\mathcal{I}(f)-\mathcal{I}_2^*(f)=O(h^4)$
for $f\in C^4(I)$\\

{\bf Example1:} $\mathcal{I}(f_1)$
\\
$$
\begin{array}{cccc}
n& E_2^*&E_2& E_3\\
 & & & \\
128&0.73(-9)&-0.55(-9)&-0.44(-8)\\
256&0;45(-10)&-0.33(-10)&-0.26(-9)\\
512&0.28(-11)&-0.21(-11)&-0.15(-10)\\
1024&0.18(-12)&-0.13(-12)&-0.95(-12)
\end{array}
$$
{\bf Example 2:} $\mathcal{I}(f_2)$
$$
\begin{array}{cccc}
n&E_2^*& E_2& E_3\\
& & & \\
128&0.14(-6)&-0.11(-6)&-0.92(-6)\\
256&0.90(-8)&-0.67(-8)&-0.52(-7)\\
512&0.56(-9)&-0.41(-9)&-0.31(-8)\\
1024&0.73(-9)&-0.52(-9)&-0.37(-8)
\end{array}
$$

{\bf (ii) QF/dQI degree 4}\\

$E_4(f)=\mathcal{I}(f)-\mathcal{I}_4(f)=O(h^6)$  for
$f\in C^6(I)$.\\

{\bf Newton-Cotes QF of degree 4:}
$E_4^*(f)=\mathcal{I}(f)-\mathcal{I}_4^*(f)=O(h^6)$
for $f\in C^6(I)$.\\

{\bf Example 1:} $\mathcal{I}(f_1)$
\\
$$
\begin{array}{ccc}
n& E_4& E_4^*\\
 & & \\
128&-0.83(-12)&1.10(-12)\\
256&-0.12(-13)&0.24(-13)\\
512&-0.18(-15)&0.37(-15)\\
1024&-0.29(-17)&0.59(-17)
\end{array}
$$
{\bf Example 2:} $\mathcal{I}(f_2)$
$$
\begin{array}{ccc}
n& E_4& E_4^*\\
 & & \\
128&0.23(-7)&-0.68(-7)\\
256&0.44(-9)&-1.04(-9)\\
512&0.73(-11)&-1.62(-11)\\
1024&0.12(-12)&-0.25(-12)
\end{array}
$$


{\bf (iii) QF/dQI degree 4:}
$E_4(f)=\mathcal{I}(f)-\mathcal{I}_4(f)=O(h^6)$  for
$f\in C^6(I)$.\\

{\bf QF/dQI degree 5:}
$E_5(f)=\mathcal{I}(f)-\mathcal{I}_5(f)=O(h^6)$ for $f\in C^6(I)$,\\

{\bf Example 1:} $\mathcal{I}(f_1)$\\
$$
\begin{array}{ccc}
n& E_4& E_5\\
 & & \\
128&-0.83(-12)&0.95(-11)\\
256&-0.12(-13)&0.14(-12)\\
512&-0.18(-15)&0.21(-14)\\
1024&-0.29(-17)&0.32(-16)
\end{array}
$$
{\bf Example 2:} $\mathcal{I}(f_2)$
$$
\begin{array}{ccc}
n& E_4& E_5\\
 & & \\
128&0.23(-7)&-0.27(-6)\\
256&0.44(-9)&-0.50(-8)\\
512&0.73(-11)&-0.83(-10)\\
1024&0.12(-12)&-0.13(-11)
\end{array}
$$

\section{Application to numerical differentiation}

Differentiating interpolation polynomials leads to classical finite
differences for the
approximate computation of derivatives. Therefore, it seems natural to
approximate
derivatives of $f$ by derivatives of $Q_df$ as long as it is possible, i.e.
up  to the order $d-1$.
The general theory will be developed elsewhere.
Here we only give results for the first derivative and $d=2,3$.
We evaluate $(Q_df)'=\sum_{j\in J} \mu_j(f) B_j'$ at points $T_n$ for $d$
even and at points $X_n$ for $d$ odd.\\

{\bf (i) Differentiation matrix for quadratics}\\

The derivation matrix $\mathcal{D}_2\in\RR^{(n+2)\times(n+2)}$ is defined
as follows: setting
$y\in\RR^{n+2}$ for the vector with components $y_j=f(t_j), j\in J$ and
$y'\in\RR^{n+2}$ for the
vector  with components $y'_j=(Q_2f)'(t_j), j\in J$, we simply write:
$$
y'=\mathcal{D}_2y
$$
$$
\mathcal{D}_2=
\begin{pmatrix}
-8/3&3&-1/3&0&0&0&\ldots&0&0\\
-7/6&11/16&13/24&-1/16&0&0&\ldots&0&0\\
1/6&-3/4&1/48&5/8&-1/16&0&\ldots&0&0\\
0&1/16&-5/8&0&5/8&-1/16&\ldots&0&0\\
\ldots&\ldots&\ldots&\ldots&\ldots&\ldots&\ldots&\ldots&\ldots\\
0&0&\ldots&1/16&-5/8&0&5/8&-1/16&0\\
0&0&\ldots&0&1/16&-5/8&-1/48&3/4&-1/6\\
0&0&\ldots&0&0&1/16&-13/24&-11/16&7/6\\
0&0&\ldots&0&0&0&1/3&-3&8/3
\end{pmatrix}
$$
\medskip

{\bf (ii) Differentiation formula for cubics}\\

The derivation matrix $\mathcal{D}_3\in\RR^{(n+1)\times(n+1)}$ is defined
as follows: setting
$y\in\RR^{n+1}$ for the vector with components $y_j=f(x_j), 0\le j\le n$
and $y'\in\RR^{n+1}$ for
the vector  with components $y'_j=(Q_3f)'(x_j),  0\le j\le n $, we obtain:
$$
y'=\mathcal{D}_3y
$$
$$
\mathcal{D}_3=
\begin{pmatrix}
-11/6&3&-3/2&1/3&0&0&\ldots&0&0\\
-1/3&-1/2&1&-1/6&0&0&\ldots&0&0\\
1/12&-2/3&0&2/3&-1/12&0&\ldots&0&0\\
0&1/12&-2/3&0&2/3&-1/12&\ldots&0&0\\
\ldots&\ldots&\ldots&\ldots&\ldots&\ldots&\ldots&\ldots&\ldots\\
0&0&\ldots&1/12&-2/3&0&2/3&-1/12&0\\
0&0&\ldots&0&1/12&-2/3&0&2/3&-1/12\\
0&0&\ldots&0&0&1/6&-1&1/2&1/3\\
0&0&\ldots&0&0&-1/3&3/2&-3&11/6
\end{pmatrix}
$$
\\
{\bf (iii) Some numerical results}\\

Again we use the two functions $\displaystyle f_1(x)=\frac{1}{1+16x^2}$ and
$f_2(x)=e^{-x}sin(5\pi
x)$ on the interval $I=[-1,1]$. For $p=1,2$, we set
$\varepsilon_p=\max_{v\in V}\vert
f_p'(v)-(Q_df_p)'(v)\vert$ where $V=T_n$ (resp. $V=X_n$) for $d=2$ (resp.
$d=3$) and
$\varepsilon_p^*=\max_{v\in V}\vert f_p'(v)-\delta f_p(v)\vert$ where
$\delta f_p(v)$ is the classical centered approximation of $f_p'(v)$ of
order $2$ (with standard
modifications at the endpoints).\\

For {\sl quadratics}, we obtain the following results
$$
\begin{array}{ccccc}
n&\varepsilon_1 & \varepsilon_1^*&\varepsilon_2&\varepsilon_2^*\\
 & & \\
64& 0.014009&0.047853&0.016143 &0.046317\\
128&0.003138&0.012079&0.003674 &0.011606\\
256&0.000767&0.003036&0.000872 &0.002902\\
512&0.000190&0.000759&0.000212 &0.00725\\
1024&0.0000475&0.0001899&0.000052&0.000181
\end{array}
$$
We see that the orders are all $O(h^2)$. However, the errors for the
derivatives of the quadratic QI
($\varepsilon_1$ and $\varepsilon_2$) are between $3$ and $4$ times less
than the errors for the
centered finite differences ($ \varepsilon_1^*$ and $\varepsilon_2^*$).\\

For {\sl cubics}, we obtain the following results

$$
\begin{array}{ccccc}
n&\varepsilon_1 & \varepsilon_1^*&\varepsilon_2&\varepsilon_2^*\\
 & & \\
64&3.0(-3) &4.7(-2)&1.0(-2)&4.7(-2)\\
128&2.0(-4)&1.2(-2) &1.4(-3)&1.2(-2)\\
256&1.3(-5)&3.0(-3)&1.8(-4)&2.9(-3)\\
512&8.0(-7)&7.6(-4)&2.4(-5)&7.2(-4)\\
1024&5.0(-8)&1.9(-4) &3.0(-6)&1.8E(-4)
\end{array}
$$
Of course, $\varepsilon_1^*$ and $\varepsilon_2^*$ are both $O(h^2)$ and
$\varepsilon_1$ and
$\varepsilon_2$ are both at least $O(h^3)$. However, for the function
$f_1$, a superconvergence
phenomenon occurs because we have $\varepsilon_1=O(h^4)$ instead of
$O(h^3)$. We shall study
this kind of results in a further paper.\\
--------------
\section{Approximating zeros of a function by those of a quadratic dQI}

Let $f$ be a continuous function defined on $I=[a,b]$. In order to locate
the zeros of $f$ in this
interval, we approximate $f$ by its $C^1$ quadratic dQI $g=Q_2f$ and we
compute the {\sl exact
zeros} of $g$: this is quite possible because $g$ is piecewise quadratic.
The complete study will
be done elsewhere. Here we take a simple example: we want to approximate
the zeros of the Legendre
polynomial $P_8(x)=$ in the interval $I=[-1,1]$. The five zeros of $P_8$
are respectively
$\{\pm x_1,\pm x_2,\pm x_3,\pm x_4\}$, with
$$
x_1=.1834346425,\; x_2=.5255324099,\; x_3=.7966664774,\; x_4=.9602898565.
$$\\
The following array gives the errors $\varepsilon_k=x_k-\bar x_k,\; 1\le
k\le 4$ where $\bar x_k$
is the zero of $g$ nearest to $x_k$.

$$
\begin{array}{ccccc}
n&\varepsilon_1 & \varepsilon_2&\varepsilon_3 & \varepsilon_4\\
 & & & &\\
16&.000543 &.003784 &.013753 &-.007841\\
32&-.000043 &.000210 &.000556 &-.001017\\
64&-.000013 &-.000012 &.000043 &.000026

\end{array}
$$
{\bf Acknowledgements:} the author thanks Professor Catterina Dagnino and
the Department of
Mathematics of the University of Turin for their kind invitation to deliver this seminar during
his stay from January 12 to 20, 2005.



{\bf Author's address: }

Paul Sablonni\`ere,\\
Centre de math\'ematiques, INSA de Rennes,\\
20 avenue des Buttes de Co\"esmes, CS 14315,\\
F-35043-RENNES C\'edex,\\
France\\

e-mail: psablonn@insa-rennes.fr


\begin{thebibliography}{99}


\bibitem{db1}
C.de Boor,
{\sl A practical guide to splines}, Revised edition. Springer-Verlag,
New-York (2001).

\bibitem{db2}
C.de Boor, Splines as linear combinations of B-splines,
{\sl Approximation Theory II}, G.G. Lorentz et al. (eds),
Academic Press, New-York (1976), 1-47.

\bibitem{ccl}
G.Chen, C.K. Chui, M.J. Lai,
Construction of real-time spline quasi-interpolation schemes,
Approx. Theory Appl. {\bf 4}\; (1988), 61-75.

\bibitem{dr}
P.J. Davis, P. Rabonowitz,
{\sl Numerical integration}. Blaisdell, Waltham (1967).

\bibitem{sv}
S.A. De Swardt, J.M. De Villiers,
Gregory type quadrature based on quadratic nodal spline interpolation.
Numer. Math. \;{\bf 85}\; (2000), 129-153.

\bibitem{dvl}
R.A. DeVore, G.G. Lorentz,
{\sl Constructive approximation}, Springer-Verlag, Berlin (1993)

\bibitem{e}
H. Engels: {\sl Numerical quadrature and cubature}. Academic Press (1980).

\bibitem{g}
W. Gautschi: Orthogonal polynomials: applications and computation. Dans {\sl
Acta Numerica 1996}, A. Iserles (ed.), CUP 1996.

\bibitem{ku}
A. Krommer, Ch. W. Ueberhuber,
{\sl Computational integration}. SIAM, Philadelphia (1998).

\bibitem{l}
V. Lampret,
An invitation to Hermite's integration and summation: a comparison between
Hermite's and
Simpson's rules.
SIAM Review {\bf 46}, No 2 (2004), 329-345.

\bibitem{lls}
B.G. Lee, T. Lyche, L.L. Schumaker,
Some examples of quasi-interpolants constructed from local spline projectors.
In {\sl Math methods for CAGD} Oslo II, 243-252.

\bibitem{ls}
T. Lyche, L.L. Schumaker,
Local spline approximation,
J. Approx. Theory {\bf 15} (1975), 294-325.

\bibitem{p}
M.J.D. Powell,
{\sl Approximation theory and methods }, Cambridge University Press (1981).

\bibitem{ps1}
P. Sablonni\`ere: Quasi-interpolantes splines sobre particiones uniformes.
First Meeting in Approximation Theory of the University of Ja\'en (Ubeda,
June 29-July 2, 2000).
Pr\'epublication IRMAR 00-38, Rennes (June 2000).

\bibitem{ps2}
P. Sablonni\`ere: On some multivariate quadratic spline quasi-interpolants
on bounded
domains. In {\sl Modern Developments in Multivariate Approximation}, W.
Haussmann et al. (eds),
 ISNM Vol. 145, Birkh\"auser Verlag (2003), 263-278.

\bibitem{ps3}
P. Sablonni\`ere: Quadratic spline quasi-interpolants on bounded domains of
${\RR}^d, d=1,2,3$.
In {\sl Spline and radial functions}, Rend. Sem. Univ. Pol. Torino, Vol.
{\bf 61} (2003), 61-78.

\bibitem{ps4}
P. Sablonni\`ere: A quadrature formula associated with a univariate
quadratic spline
quasi-interpolant. Pr\'epublication IRMAR, Rennes, April 2005 (submitted).

\bibitem{ssb}
P. Sablonni\`ere, D. Sbibih: Integral spline operators exact on
polynomials. Approx. Theory
Appl. 10, No 3 (1994), 56-73.

\bibitem{sc}
L.L. Schumaker,
{\sl Spline functions: basic theory}, John Wiley and Sons, New-York (1981).

\end{thebibliography}
\end{document}